 \documentclass[12pt]{article}
\usepackage{mathrsfs}
\usepackage{amsfonts}
\usepackage{latexsym}

\usepackage[dvips]{graphicx}
\usepackage{geometry,amsthm,graphics,amssymb,amsmath,enumerate,latexsym,tabularx,shapepar}
\usepackage[all,2cell,dvips]{xy} \UseAllTwocells \SilentMatrices
\date{}

\begin{document}
\title{Subalgebras and finitistic dimensions of Artin algebras$^\star$ }
\author{{\small Aiping Zhang$^{1}$,   Shunhua Zhang $^{2}$}\\
{\small 1 School of Mathematics and Statistics,\ Shandong University
at Weihai,\ 264209,\ China }\\{\small 2 School of Mathematics,\
Shandong University,\ Jinan
250100,\ China }\\
}

\pagenumbering{arabic}

\maketitle

\footnote{ {\it Email addresses}: pingping326@163.com(A.Zhang ),
shzhang@sdu.edu.cn(S.Zhang)}

\footnote{ $^\star$Supported by the NSF of China (Grant No.
10771112) and NSF of Shandong Province (Grant No. Y2008A05).}

\begin{center}
 \begin{minipage}{120mm}
   \small\rm
   {\bf  Abstract}\ \ Let $A$  be an Artin algebra. We investigate
   subalgebras of $A$ with certain conditions and obtain some
classes of algebras whose  finitistic dimensions are finite.
\end{minipage}
\end{center}
\begin{center}
  \begin{minipage}{120mm}
   \small\rm
   {\it Keywords:}{ \ \ Artin algebras,\ representation dimension,\ finitistic
   dimension}
\end{minipage}
\end{center}
\vskip -0.3cm
\begin{center}
\begin{minipage}{120mm}
   \small\rm
   {\bf MSC(2000):}\ {\small\bf16G10,  18G20}
\end{minipage}
\end{center}

\section {Introduction}

\vskip 0.2in

Let $A$  be  an  Artin algebra, $A$-mod the category of finitely
generated left $A$-modules, and $A$-ind a full subcategory of
$A$-mod containing exactly one representative of each isomorphism
class of indecomposable $A$-modules. We denote the projective
dimension of an $A$-module $X$ by ${\rm pd}_A\ X$.

\vskip 0.2in

 Let $A$  be  an  Artin algebra. Recall from [1] that the finitistic dimension  of
 $A$, denoted by ${\rm fin.dim}\ A$, is defined as
$$
{\rm fin.dim}\ A = {\rm sup}\{{\rm pd}_{A}\ M \ |\ M\in A-{\rm
mod},\ {\rm pd}_{A}\ M<\infty\}.
$$
The finitistic dimension conjecture claims that every Artin algebra
has a finite finitistic dimension.

\vskip 0.2in

So far, only a few classes of algebras were known to have finite
finitistic dimensions. For example, monomial algebras [3], algebras
where the cube of the radical is zero [4], and the algebras given in
[5,6,7,8,9,10,11,12]. However, the finitistic dimension conjecture
is still open and it is far from to be proven.

\vskip 0.2in

Let $A$ be an Artin algebra, and $0\rightarrow \ _AA\rightarrow
I_0\rightarrow I_1\rightarrow \cdots$ be the minimal injective
resolution of $A$. Nakayama conjectured in [13] that $A$ is a
self-injective algebra whenever all $I_j$ is projective. Up to now,
Nakayama conjecture is still open. It is well known that finitistic
dimension conjecture implies Nakayama conjecture, and this motivated
further research on finitistic dimension conjecture. We refer to
[8,9,10,11,12,14] for the background and some new progress about
this conjecture.

\vskip 0.2in

According to Xi in [9], the finitistic dimension conjecture is
equivalent to the following statement: if $B$ is a subalgebra of $A$
such that rad $B$ is a left ideal in $A$,  then $B$ has finite
finitistic dimension whenever $A$ has finite finitistic dimension.

\vskip 0.2in

 In this paper, we investigate the finitistic dimensions
 of Artin algebras by using Igusa-Todorov function defined in [5] and
 obtain some classes of algebras with finite finitistic dimensions.
 The paper is arranged as follows. In Section 2 we collect some definitions
 and results needed for our research, and give a different proof for
 a well known fact (Theorem 2.5).  In Section 3, we  obtain some classes of
 algebras with finite finitistic dimensions, which gives a partly positive answer
 to the question 2 mentioned in [10].

\vskip 0.2in

\section {Preliminaries}

\vskip 0.2in

Throughout this paper, we always assume that $A$ is an Artin
algebra. We denote the global dimension of $A$ by gl.dim $A$ and the
Jacobson radical of $A$ by rad $A$. For an $A$-module $M$, we denote
by add $M$ the full subcategory having as objects the direct sums of
indecomposable summands of $M$,  by $\Omega^i M$ the $i$th syzygy of
$M$.  Then, $\cal  P$ = add $_A$$A$ is the full subcategory
consisting of all finitely generated projective A-modules, and $\cal
I$ = add $_A$$DA$ is the full subcategory consisting of all finitely
generated injective A-modules, where $D: A-{\rm mod}\rightarrow
A^{op}-{\rm mod}$ is the standard duality, and $A$$^{op}$ is the
opposite algebra of $A$. Given two homomorphisms $f:
L$$\longrightarrow$$M$ and $g: M$$\longrightarrow$$N$, the
composition of $f$ and $g$ is denoted by $gf$. We follow the
standard terminology and notation used in the representation theory
of algebras, see [15] and [16].

\vskip 0.2in

An $A$-module $V$ is called a $generator-cogenerator$ if every
indecomposable projective module and every indecomposable injective
module is isomorphic to a summand of $V$. Recall from [18] that the
number
$$
{\rm rep.dim} \ A= {\rm inf} \ \{ \ {\rm gl.dim\ End}_A(V) \ \ | \ \
V \ {\rm is \  a  \ generator-cogenerator} \ \}
$$
is called the $representation\ dimension$ of an Artin algebra $A$.

 \vskip 0.2in

{\bf Lemma 2.1.}\ \ {\it Let $M$ be an $A$-module and  there is an
exact sequence $0\rightarrow X_{s}\rightarrow\cdots \rightarrow
X_{1}\rightarrow X_{0}\rightarrow M\rightarrow 0$. If \ ${\rm pd}_A\
X_i \leq k$, $i=  0,\cdots, s$, then \  ${\rm pd}_A\ M \leq s+k$.}
\hfill$\Box$

\vskip 0.2in

The following two Lemmas proved in [18] and [9] will be used later.

\vskip 0.2in

{\bf Lemma 2.2.}  {\it Let $V$ be a generator-cogenerator of
$A$-{\rm mod} and $n \geq 3$  an integer. The following two
statements are equivalent:

\vskip 1mm

 {\rm (1)}  For any $X\in A$-{\rm ind}, there is an exact sequence
$$
0 \rightarrow V_{n-2}\rightarrow\cdots\rightarrow V_{1}\rightarrow
V_{0} \rightarrow X \rightarrow 0
$$
with $V_i\in {\rm add}\ (_AV)$ for $j =  0,\cdots, n-2$, such that
$$
0\rightarrow {\rm Hom}_{A}(V,V_{n-2})\rightarrow\cdots\rightarrow
{\rm Hom}_{A}(V,V_{1}) \rightarrow {\rm Hom}_{A}(V,V_{0})
\rightarrow {\rm Hom}_{A}(V,X) \rightarrow 0
$$
is exact.

\vskip 1mm

 {\rm (2)}  ${\rm gl.dim\ End}\ _AV \leq  n$. }               \hfill$\Box$

\vskip 0.2in

{\bf Lemma 2.3.}  {\it Suppose $B$ is a subalgebra of $A$ such that
\ ${\rm rad}\ B$  is a left ideal in $A$. For any $B$-module $X$ and
integer $i \geq 2$, there is a projective $A$-module $Q$ and an
$A$-module $Z$ such that \ $\Omega_{B}^{i}(X) \cong\Omega_{A}(Z)
\oplus Q$ as $A$-modules. If \ ${\rm rad}\ B$ is an ideal in $A$,
then there is an exact sequence of $A$-modules
$$
0\longrightarrow \Omega_B^i(X) \longrightarrow
\Omega_A^2(Y) \oplus P \longrightarrow S\longrightarrow 0,
$$
where $P$ is projective, and $S$ is an $A$-module such that $_BS$ is
semisimple. In particular, if \ ${\rm rad}\ B={\rm rad} A$, the
module $S$ is even a semisimple $A$-module.} \hfill$\Box$

\vskip 0.2in

Let $K(A)$ be the free abelian group with the basis of
non-isomorphism classes of non-projective indecomposable $A$-modules
in $A$-mod. Igusa and Todorov in [5] define a function $\psi_A$ on
$K(A)$, which depends on the algebra $A$ and take values of
non-negative integers. Now $\psi_A$ is called Igusa-Todorov
function, and it is a powerful tool to show the finiteness of the
finitistic dimensions. The following lemma collects some important
properties of this Igusa-Todorov function.

\vskip 0.2in

{\bf Lemma 2.4.}  {\it Let $A$ be an Artin algebra, and $\psi_A$ be
the corresponding Igusa-Todorov function. Let $M,X,Y,Z$ be
$A$-modules in $A$-mod.}

\vskip 0.1in

(1)\ {\it $\psi_A (M)={\rm pd}\ M$ provided ${\rm pd}\ M<\infty$. }

\vskip 0.1in

(2)\ {\it $\psi_A (X)\leq \psi_A (X\oplus Y)$.}

\vskip 0.1in

(3)\ {\it If $0\rightarrow X\rightarrow Y\rightarrow Z\rightarrow 0$
is an exact sequence in $A$-mod and ${\rm pd}\ Z<\infty$, then ${\rm
pd}\ Z\leq \psi_A (X\oplus Y)+1.$}

\vskip 0.2in

Let $\cal X$ be a full subcategory of $A$-mod. When we say that
$\cal X$ is a full subcategory, we always mean that $\cal X$ is
closed under direct summands. We denote by gen $\cal X$ (cogen $\cal
X$) the full subcategory of $A$-mod generated (cogenerated) by $\cal
X$, see [17] and [16]. If $\cal X$ = $\{M\}$, we set ${\cal X}=M$
and denote gen $\cal X$ (cogen $\cal X$) by gen $M$ (cogen $M$). If
$\cal X$ contains only finite non-isomorphic indecomposable
$A$-modules, we call $\cal X$ is of finite type.

\vskip 0.2in

It has been shown in [20] that ${\rm rep.dim} \ A$ is at most $3$
whenever ${\rm gen}\ DA$ is finite, then in this case, according to
[5], ${\rm fin.dim}\ A$ is finite. Now, we give a different proof
for this result by using Igusa-Todorov function.

\vskip 0.2in

{\bf Theorem 2.5.}  {\it Let $A$ be an Artin algebra. If ${\rm gen}\
DA$ is of finite type, then ${\rm fin.dim}\ A$ is finite.}

\vskip 0.1in

{\bf Proof.}  Let $X$ be an $A$-module with finite projective
dimension. Let $i: X \longrightarrow  E(X)$ be the injective
envelope of $X$, then we have an exact sequence 0$\longrightarrow
X\stackrel{i}\longrightarrow E(X)\longrightarrow {\rm coker} \
i\longrightarrow 0$.

 We may assume that
$M_{1}$ , $\cdots$, \ $M_{t}$ are a complete list of pairwise
non-isomorphic indecomposable $A$-modules in gen $DA$. Since the
modules $E, {\rm coker}  \   i$ lie in {\rm gen} $DA$, we may write
$E(X) = \bigoplus\limits_{i=1}^{t}M_{i}^{t_i}$,  ${\rm coker}\
i=\bigoplus\limits_{j=1}^{t}M_{j}^{s_j}$.  We denote by $a= {\rm
max}\{\  t_i+s_i \  \}$.  By Lemma 2.4, we know that
$$\begin{array}{ccl}
  {\rm pd}\ X&\leq & \psi_A(\Omega(E(X)\oplus {\rm coker}\ i))+1\\
     &= & \psi_A(\Omega (E(X))\oplus\Omega({\rm coker}\ i))+1\\
     &= & \psi_A(\Omega(\bigoplus\limits_{i=1}^{t}M_{i}^{t_i})
     \oplus\Omega(\bigoplus\limits_{j=1}^{t}M_{j}^{s_j}))+1\\
     &= & \psi_A(\bigoplus\limits_{i=1}^{t}\Omega(M_{i})^{t_i}
     \oplus\bigoplus\limits_{j=1}^{t}\Omega(M_{j})^{s_j})+1\\
     &= & \psi_A(\bigoplus\limits_{i=1}^{t}\Omega(M_{i})^{t_i+s_i})+1\\
     &\leq & \psi_A(\bigoplus\limits_{i=1}^{t}\Omega(M_{i})^{a})+1\\
     &= & \psi_A(\bigoplus\limits_{i=1}^{t}\Omega(M_{i}))+1.
     \end{array}
     $$

Thus
 fin.dim $A\leq\psi_A(\bigoplus\limits_{i=1}^{t}\Omega(M_{i}))+1$,
 it follows that the finitistic dimension
of $A$ is finite.     \hfill$\Box$

\vskip 0.2in

\section {Main results }

 \vskip 0.2in

In this section, we investigate the finitistic dimensions of
subalgebras of an Artin algebra  with certain conditions and give
some examples to show how our results are applied.

\vskip 0.2in

 We replace the condition gl.dim $A$ $\leq$ 1 of Theorem 3.3
in [10] by rep.dim $A$ $\leq$ 3 and obtain the following result.

\vskip 0.2in

{\bf Theorem 3.1.}\ \  {\it Let $A_{0} = B\subseteq A_{1}\subseteq
\cdots \subseteq  A_{s-1} \subseteq  A_{s} = A$ be a chain of
subalgebras of $A$,  ${\rm rad}\ (A_{i-1})$ is a left ideal in
$A_{i}$ for all $i$ and  ${\rm pd}_{A_{i-1}}\ A_{i} < \infty$ for
all $1\leq  i  \leq  s-1$. If \ ${\rm rep.dim}\ A \leq 3$, then
${\rm fin.dim}\ B < \infty$.}

\vskip 1mm

{\bf Proof.}\ \ According to the proof of Theorem 3.1 in [10]. We
know that pd$_{B}\ A_{j}< \infty$ for all $1 \leq j \leq s-1$.
Suppose $M$ is a $B$-module, pd $_B M < \infty$. We denote by
$\Omega_{i}$ the first syzygy operator of $A_{i}$-modules. By Lemma
2.3, $\Omega_0^2(M)$ is an $A_{1}$-module. Similarly
$\Omega_j^2\cdots\Omega_1^2\Omega_0^2(M)$ is an $A_{j+1}$-module, we
have the following exact sequences:

$0\rightarrow\Omega_0^2(M)\rightarrow P_{0}(1)\rightarrow
P_{0}(0)\rightarrow M\rightarrow 0,$

\vskip 1mm

$0\rightarrow\Omega_1^2\Omega_0^2(M)\rightarrow
P_{1}(1)\rightarrow P_{1}(0)\rightarrow\Omega_0^2(M) \rightarrow
0$,

\vskip 1mm

$0\rightarrow\Omega_2^2\Omega_1^2\Omega_0^2(M)\rightarrow
P_{2}(1)\rightarrow P_{2}(0)
\rightarrow\Omega_1^2\Omega_0^2(M)\rightarrow 0$,

$\vdots$

$0\rightarrow\Omega_{s-1}^2\Omega_{s-2}^2\cdots\Omega_0^2(M)\rightarrow
P_{s-1}(1) \rightarrow
P_{s-1}(0)\rightarrow\Omega_{s-2}^2\cdots\Omega_0^2(M)\rightarrow
0,   \ \ \ (**) $

\vskip 1mm

where $P_0(1),\ P_0(0)$ are projective $A_0$-modules, $P_1(1),\
P_1(0)$ are projective $A_1$-modules, $P_2(1),\ P_2(0)$ are
projective $A_2$-modules , $\cdots, P_{s-1}(1),\ P_{s-1}(0)$ are
projective $A_{s-1}$-modules.

Thus we have the following long exact sequence

 $(1)\ \
0\rightarrow\Omega_{s-1}^2\Omega_{s-2}^2\cdots\Omega_0^2(M)\rightarrow
P_{s-1}(1) \rightarrow P_{s-1}(0)\rightarrow\cdots\rightarrow
P_{1}(1)\rightarrow P_{1}(0) \rightarrow P_{0}(1)\rightarrow
P_{0}(0)\rightarrow M\rightarrow 0$.

 By ($**$), we know that $B$-module
 $\Omega_{s-1}^2\Omega_{s-2}^2\cdots\Omega_0^2(M)$ has finite
projective dimension.

It follows from Lemma 2.2 and the inequality
 rep.dim $A \leq 3$  that there exists a
generator-cogenerator $V$ for $A-$mod, such that for any
$A$-module $X$, there is an exact sequence $0 \longrightarrow V_1
\longrightarrow V_0 \longrightarrow X \longrightarrow 0,$ with
$V_1, V_0 \in$ add $V$, such that
$$
0\rightarrow {\rm Hom}_{A}(V,V_1)\rightarrow {\rm Hom}_{A}(V,V_0)
\rightarrow {\rm Hom}_{A}(V,X) \rightarrow 0
$$
is exact.

Obviously $\Omega_{s-1}^2\Omega_{s-2}^2\cdots\Omega_0^2(M)$ is an
$A$-module and there is a short exact sequence
\begin{center}
$0\longrightarrow V_{1}\longrightarrow V_{0}\longrightarrow
\Omega_{s-1}^2\Omega_{s-2}^2\cdots\Omega_0^2(M)\longrightarrow 0$
\end{center}
of $A$-module with $V_{1}$, $V_{0}$ $\in$ add $V$. By Lemma 2.4, we
know that
$$\begin{array}{ccl}
 {\rm pd} _B\Omega_{s-1}^2\Omega_{s-2}^2\cdots\Omega_0^2(M)&\leq & \psi_B(V_1\oplus V_0)+1\\
     &\leq& \psi_B(V)+1.\\
     \end{array}
$$
By the long exact sequence (1), we have
$$\begin{array}{ccl}
 {\rm pd} _BM&\leq & 2s+{\rm max}\{\ {\rm
 pd}_B\Omega_{s-1}^2\Omega_{s-2}^2\cdots\Omega_0^2(M),
 {\rm pd}_B(P_j(i))|^{i=0,\ 1}_{ j=1, \cdots,\ s-1} \ \}\\
     &\leq& 2s+{\rm max}\{ \ \psi_{B}(V)+1,\ {\rm
     pd}_B(A_j)|j=1,\ \cdots,\ s-1 \ \},
    \end{array}
     $$
Thus fin.dim $B$ is finite.
    \hfill$\Box$

 \vskip 0.2in

When s = 2, we obtain the following consequence.

 \vskip 0.2in

{\bf Corollary 3.2.}  {\it Let $C \subseteq  B \subseteq A$ be a
chain of subalgebras of an Artin algebra $A$ such that \ ${\rm rad}\
C$\ is a left ideal in $B$, ${\rm rad}\ B$ is a left ideal in $A$.
If ${\rm pd}_C\ B < \infty$ and ${\rm rep.dim}\ A \leq 3$, then
${\rm fin.dim}\ C$ is finite.}  \hfill$\Box$

 \vskip 0.2in

{\bf Example 1.}\ \ Let $A$ be an algebra (over a field ) given by
the following quiver with relations: $cd=ef$, $a^4=ba=0$.

\begin{figure}[htbp]
\begin{center}
\includegraphics[width=8cm,height=3cm]{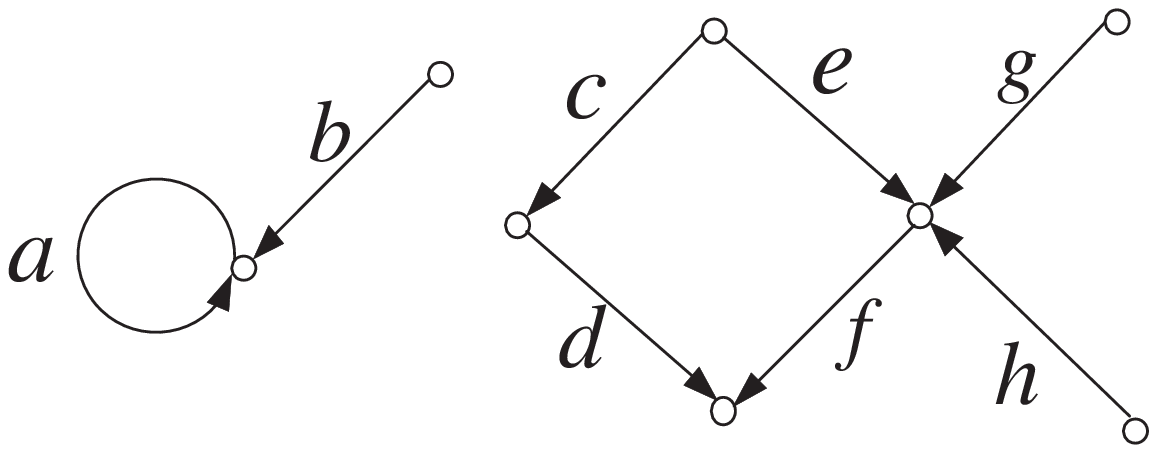}
\end{center}
\end{figure}

Now we use the method "gluing of idempotents" to construct
subalgebras of $A$, which have the same radical rad $A$, see [10]
for details.

Let $B$ be the subalgebra of $A$ given by the following quiver with
relations: $cd=ef$, $a^4=ba=ca=bd=ad=0$:

\begin{figure}[htbp]
\begin{center}
\includegraphics[width=8cm,height=3cm]{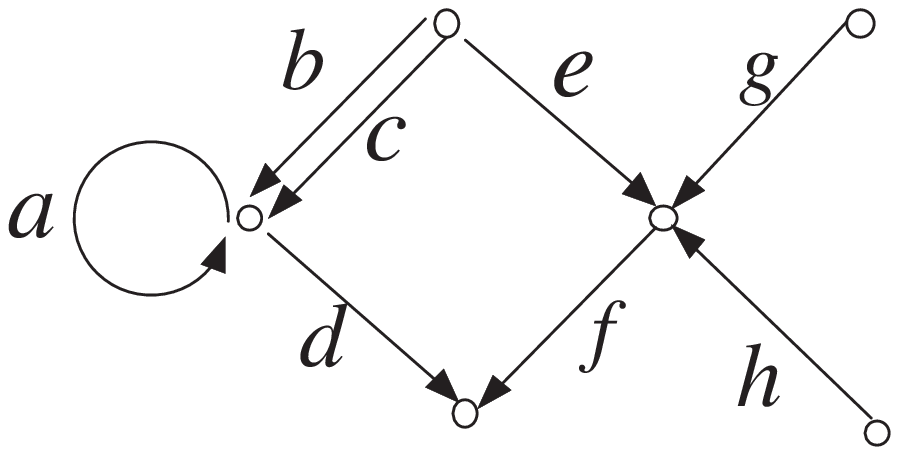}
\end{center}
\end{figure}

\vskip 0.8in

Now we consider the subalgebra $C$ of $B$, which is given by quiver
and relations $cd=ef$, $a^4=ba=ca=bd=ad=0$, $gf=hf$.

\begin{figure}[htbp]
\begin{center}
\includegraphics[width=8cm,height=3cm]{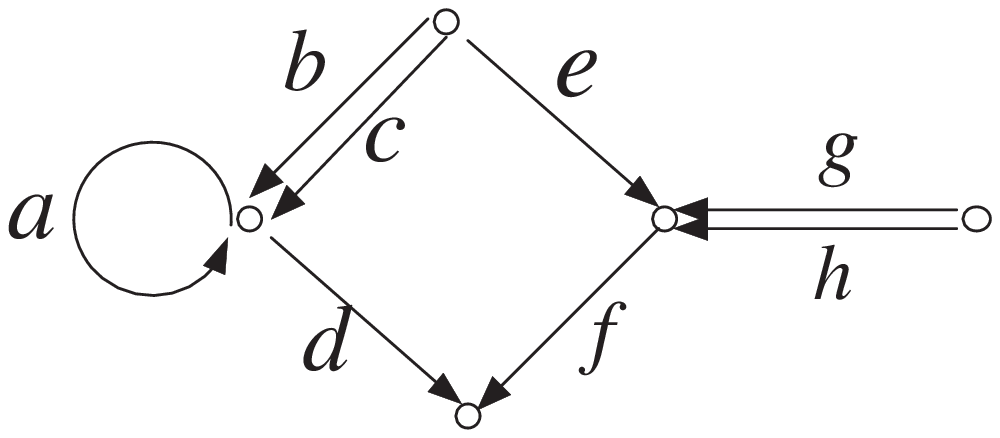}
\end{center}
\end{figure}

By Corollary 2.4 in [19],  rep. dim $A \leq 3$, and it is easy to
see that ${\rm pd}_C B=2$. Then we have that fin. dim $C< \infty$ by
Corollary 3.2.

 \vskip 0.2in

{\bf Theorem 3.3.}  {\it Let $B$ be a subalgebra of an Artin algebra
$A$ such that ${\rm rad}\ B$ is an ideal in $A$. If ${\rm add} \{
\Omega_A^2(M)\ |\ M \in  A-{\rm mod}\}$ is of finite type, then
${\rm fin.dim}\ B$ is finite.}

\vskip 1mm

{\bf Proof.}  Suppose $X$ is a $B$-module with finite projective
dimension, by Lemma 2.3, there is an exact sequence of $A$-modules:
\begin{center}
$(*)\ \ \ 0\longrightarrow \Omega_B^2(X)\longrightarrow
\Omega_A^2(Y)\oplus P \longrightarrow S\longrightarrow 0,$
\end{center}
where $_A$$P$ is projective, and $S$ is an $A$-module such that
$_B$$S$ is semisimple.

Since add$\{$$\Omega_A^2$($M$)$\mid$$M$$\in $$A$-mod$\}$ is of
finite type, we may assume that $M_{1}$,\ $\cdots$,\ $M_{t}$ are a
complete list of pairwise non-isomorphic indecomposable $A$-modules.
Obviously $\Omega_A^2(Y)$ lie in add
$\{$$\Omega_A^2$($M$)$\mid$$M$$\in $$A$-mod$\}$, so we may write
$\Omega_A^2(Y)=\bigoplus\limits_{i=1}^{t}M_i^{t_i}$. By Lemma 2.4
and ($*$), we know that
$$\begin{array}{ccl}
{\rm pd}\ _B X&\leq & {\rm pd}\ \Omega_B^2(X)+2\\
     & \leq & \psi_B(\Omega_B(\Omega_A^2(Y)\oplus P \oplus S))+3\\
     &= & \psi_B((\Omega_B(\bigoplus\limits_{i=1}^{t}M_{i}^{t_i})
     \oplus\Omega_B(P)\oplus\Omega_B(S))+3\\
     &\leq & \psi_B(\Omega_B(M_1)\oplus\cdots\oplus\Omega_B(M_t)
     \oplus\Omega_B(A)\oplus\Omega_B(B/{\rm rad}B))+3.
\end{array}
$$
Thus fin.dim $B$ is finite. \hfill$\Box$

\vskip 0.2in

{\bf Corollary 3.4.}  {\it Let $B$ be a subalgebra of an Artin
algebra $A$ such that ${\rm rad} B$ is an ideal in $A$. If ${\rm
gl.dim} A \leq  2$, then ${\rm fin.dim} B$ is finite.} \hfill$\Box$

 \vskip 0.2in

{\bf Theorem 3.5.}  {\it Let $C \subseteq B \subseteq A$ be a chain
of subalgebras of an Artin algebra $A$ such that ${\rm rad}\ C$ is a
left ideal in $B$, ${\rm rad}\ B$ is a left ideal in $A$. If ${\rm
cogen}\ A$ is of finite type, then ${\rm fin.dim}\ C$ is finite.}

{\bf Proof.}\ \ Suppose $X$ is a $C$-module with finite projective
dimension. By Lemma 2.3, $\Omega_{C}^2(X)$ is a $B$-module, we have
the following exact sequence
$$0\longrightarrow \Omega_B\Omega_{C}^2(X)\longrightarrow
P\longrightarrow \Omega_{C}^2(X)\longrightarrow 0,
$$ of $B$-modules
where $P$ is a projective $B$-module. By Lemma 2.3, there is a
$B$-module $Y$ and a projective $B$-module $Q^{\prime}$ such that
$\Omega_{C}^2(X)=\Omega_B(Y)\oplus Q^{\prime}$. Thus the above exact
sequence can be rewritten as $$0\longrightarrow
\Omega_{B}^2(Y)\longrightarrow P\longrightarrow
\Omega_{C}^2(X)\longrightarrow 0,
$$  by Lemma 2.3 again, there is an
$A$-module $Z$ and a projective  $A$-module $Q$ such that
$\Omega_{B}^2$(Y)=$\Omega$$_A$(Z)$\oplus$ $Q$, so we have the
following exact sequence
$$0\longrightarrow \Omega_A(Z)\oplus Q\longrightarrow
P\longrightarrow \Omega_{C}^2(X)\longrightarrow 0.$$Since cogen $A$
is of finite type, we may assume that $M_{1}$,\ $\cdots$,\ $M_{t}$
are a complete list of pairwise non-isomorphic indecomposable
$A$-modules, $\Omega_A(Z) \in {\rm cogen}\ A$, $Q\in {\rm cogen}\
A$,  we may write
$\Omega_A(Z)=\bigoplus\limits_{i=1}^{t}M_{i}^{t_i}$,
$Q=\bigoplus\limits_{j=1}^{t}M_{j}^{s_j}$, By Lemma 2.4, we know
that
$$\begin{array}{ccl}
{\rm pd}\ _{C}X&\leq & {\rm pd}\ \Omega_{C}^2(X)+2\\
     &\leq & \psi_C(\Omega_A(Z)\oplus Q\oplus P)+3\\
     &\leq & \psi_C(\bigoplus\limits_{i=1}^{t}M_{i}^{t_i}\oplus\bigoplus\limits_{j=1}^{t}M_{j}^{s_j}\oplus B)+3\\
     &\leq & \psi_C(M_1\oplus\cdots\oplus M_t\oplus B)+3.
      \end{array}
     $$
Thus ${\rm fin.dim}\ C$ is finite.\hfill$\Box$

\vskip 0.2in

{\bf Remark.}\   we should point that Corollary 3.2 and Theorem 3.5
give partial answer to the question 2 in [10].

\vskip 0.2in

 {\bf Corollary 3.6.}\ \ {\it Let $C\subseteq B\subseteq A$
 be a chain of subalgebras of an Artin algebra $A$ such
that ${\rm rad}\ C$ is a left ideal in $B$, ${\rm rad}\ B$ is a left
ideal in $A$. If $A$ is a hereditary Artin algebra, then ${\rm
fin.dim}\ C$ is finite.}

\newpage

{\bf Example 2.}\ \ Let $A$ be an algebra (over a field ) given by
the following quivers with relations: $a^4=0$, $cd=ef$.

\begin{figure}[htbp]
\begin{center}
\includegraphics[width=8cm,height=3cm]{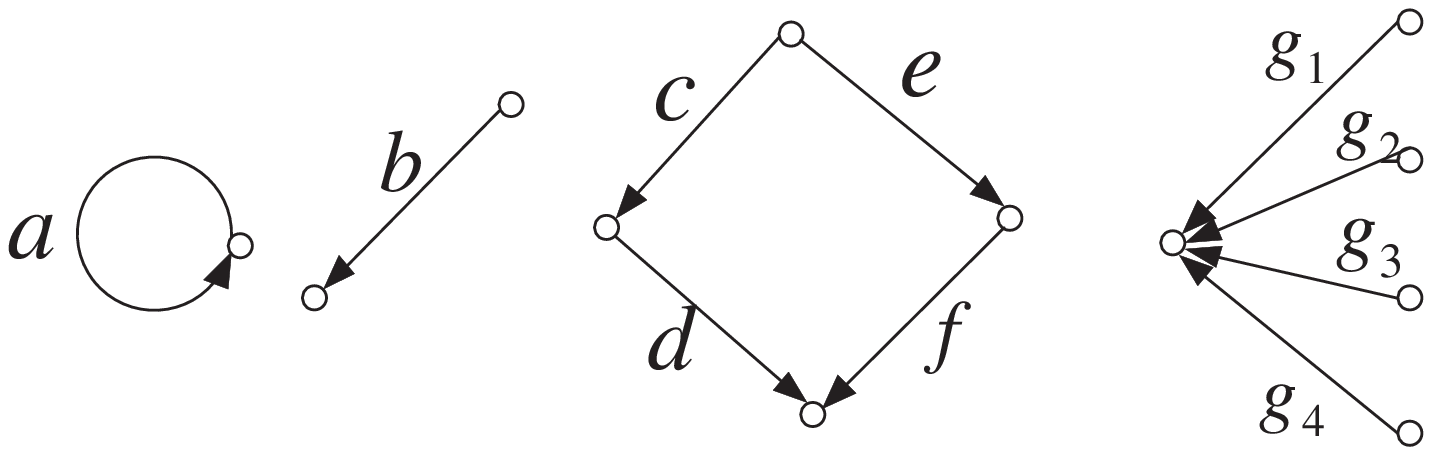}
\end{center}
\end{figure}

We use the method "gluing of idempotents" (see [10]) to construct
subalgebras of $A$ as following.

\vskip 0.1in

Let $B$ be the subalgebra of $A$ given by the following quivers with
relations: $a^4=ba=0$, $ca=ef$, $g_1f=g_2f$, $g_3f=g_4f$.

\begin{figure}[htbp]
\begin{center}
\includegraphics[width=8cm,height=3cm]{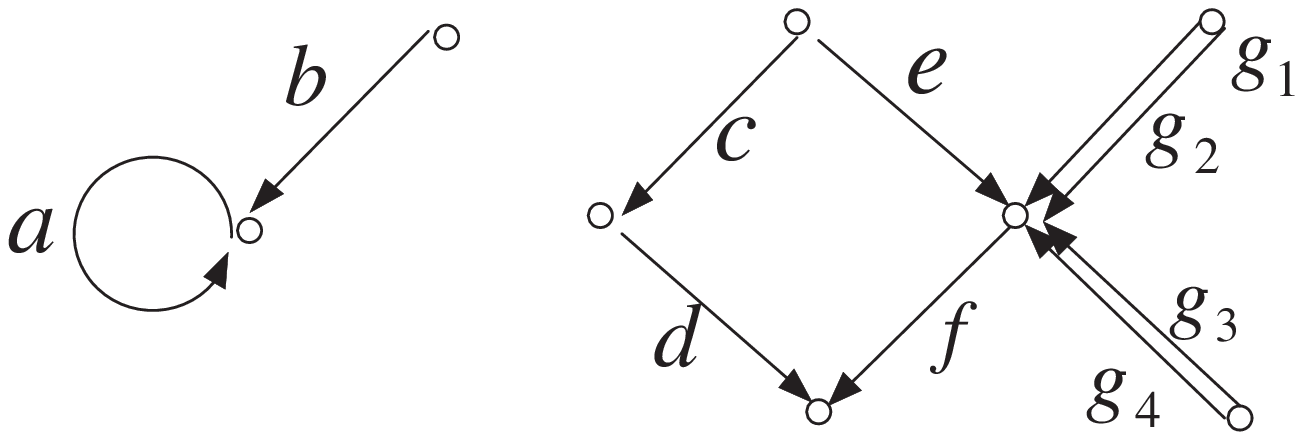}
\end{center}
\end{figure}

Let $C$ be the subalgebra of $B$ given by the following quiver with
relations $a^4=ba=bd=ca=ad=0$, $cd=ef$, $g_1f=g_2f=g_3f=g_4f$.

\begin{figure}[htbp]
\begin{center}
\includegraphics[width=8cm,height=3cm]{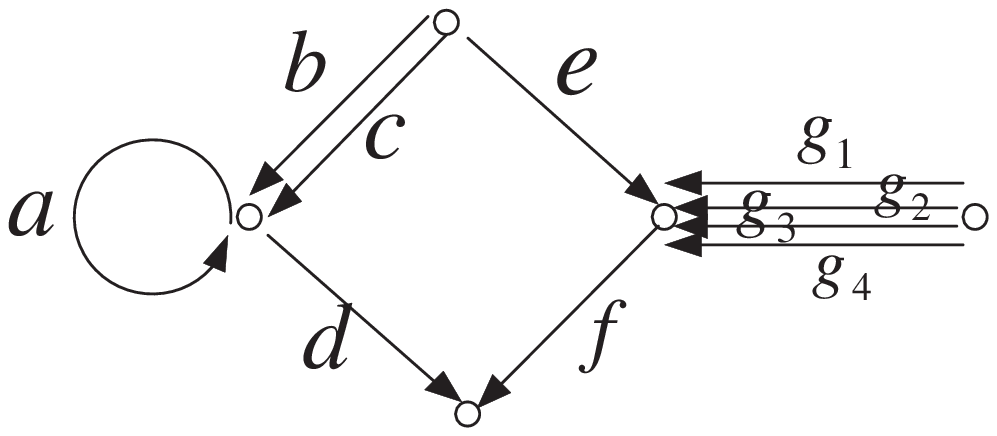}
\end{center}
\end{figure}

Then we get a chain of subalgebras of $A$, $C \subseteq B \subseteq
A$,  such that ${\rm rad}\ C$ is a left ideal of $B$, ${\rm rad}\ B$
is a left ideal of $A$. It is easy to see that ${\rm cogen}\ A$ is
of finite type. By Theorem 3.5, we know that ${\rm fin.dim}\ C <
\infty$.

\vskip 0.9in

 {\bf Acknowledgements}\ \ The authors would like to thank the referee for his or her
 helpful suggestions. They also want to thank Wenxu Ge and Hongbo Lv
 for fruitful discussions. In particular, the Example 1 is pointed
to us by Hongbo Lv and the quivers in Example 1 and Example 2 are
edited by Wenxu Ge.

\begin{description}

\item{1}\ Bass H. Finitistic dimension and a homological
generalization of semiprimary rings. Trans.Amer.Math.Soc, 1960,95:
466-488

\item{2}\ Zimmermann-Huisgen B. Homological domino effects and
the first finitistic dimension conjecture. Invent.Math, 1992,108:
369-383

\item{3}\ Green E, Kirkman E, Kuzmanovich J. Finitistic
dimension of finite dimensional monomial algebras. J.Algebra.
1991,136: 37-51

\item{4}\ Green E, Zimmermann-Huisgen B. Finitistic dimension
of artinian rings with vanishing radical cube. Math.Z. 1991,206:
505-526

\item{5}\ Igusa k, Todorov G. On the finitistic global
dimension conjecture for artin algebras, In:Representations of
algebras and related topics,201-204. Fields Inst. Commun., 2002,45,
Amer.math.Soc., Providence, RI, 2005.

\item{6}\ Auslander M, Reiten I. Applications of
contravariantly finite subcategories. Adv.Math. 1991,86: 111-152

\item{7}\ Wang Y. A note on the finitistic dimension conjecture.
Comm.Algebra.  1994, 22(7): 2525-2528

\item{8}\ Xi C C.  Representation dimension and quasi-hereditary
algebras. Adv.Math. 2002, 168:193-212

\item{9}\ Xi C C, On the finitistic dimension conjecture I: \
related to representation-finite algebras, J.Pure Appl.Algebra.
2004,193: 287-305

\item{10}\ Xi C C, On the finitistic dimension conjecture II: \
related to finite global dimension, Adv.Math. 2006,201: 116-142

\item{11}\ Xi C C, On the finitistic dimension conjecture III:\
related to the pair $eAe \subseteq A$,  J.Algebra,
2008,319(9):3666-3688

\item{12}\ Zhang A P, Zhang S H. On the finitistic  dimemsion
conjecture of Artin algebras. J.Algebra, 2008,320:253-258

\item{13}\ Nakayama T, On algebras with complete homology,  Abh.
Math. Sem. Univ. Hamburg. 1958,22: 300-307.

\item{14}\ Zimmermann-Huisgen B. The finistic dimemsion
conjecture--a tale of 3.5 decates. Abelian group and modules
(Padova, 1994), Math.Appl., Vol.343, Kluwer Academic Publishers,
Dordrecht, 1995, 501-517.

\item{15}\ Auslander M, Reiten I,Smal$\phi$  S.O. Representation
Theory of Artin Algebras. Cambridge University Press,\ Cambridge,
1995.

\item{16}\ Assem I, Simson D, Skowronski A. Elements of the
representation theory of associative algebras, Vol. 1, Cambridge
Univ. Press, 2006.

\item{17} \ Anderson F W, Fuller K R. \ Rings and categories of
modules. Springer 1973.

\item{18}\ Auslander M. Representation Dimension of Artin
Algebras. Queen Mary College Mathematics Notes, Queen Mary College,
London, 1971.

\item{19}\ Coelho F.U.,Platzeck M.I.  On the representation  dimemsion
 of some classes of  algebras. J.Algebra, 2004, 275(2):615-628.

\item{20}\ Ringel C.M. The torsionless modules of an Artin algebra,
preprint 2008.

\end{description}

\end{document}